\documentclass[12pt]{article}%
\usepackage{amsmath}
\usepackage{amsfonts}
\usepackage{amssymb}
\usepackage{graphicx}%
\newtheorem{theorem}{Theorem}

\newtheorem{conclusion}[theorem]{Conclusion}

\newtheorem{lemma}[theorem]{Lemma}

\newtheorem{remark}[theorem]{Remark}

\begin{document}

\title{Asymptotic analysis of the Askey-scheme I: from Krawtchouk to Charlier}
\author{Diego Dominici \thanks{e-mail: dominicd@newpaltz.edu}\\Department of Mathematics\\State University of New York at New Paltz\\75 S. Manheim Blvd. Suite 9\\New Paltz, NY 12561-2443\\USA\\Phone: (845) 257-2607\\
Fax: (845) 257-3571
}
\maketitle

\newpage

\begin{abstract}
We analyze the Charlier polynomials $C_{n}(x)$ and their zeros asymptotically
as $n\rightarrow\infty.$ We obtain asymptotic approximations, using the limit
relation between the Krawtchouk and Charlier polynomials, involving some
special functions. We give numerical examples showing the accuracy of our formulas.

\end{abstract}


MSC-class: 33C45 (Primary) 34E05, 33C10 (Secondary)

\newpage

\section{Introduction}

The Charlier polynomials $C_{n}(x)$ \cite{charlier} are defined by
\begin{equation}
C_{n}(x)=\,_{2}F_{0}\left(  \left.
\begin{array}
[c]{c}%
-n,-x\\
-
\end{array}
\right\vert -\frac{1}{a}\right)  \label{Cndef}%
\end{equation}
where $x\geq0,\ n=0,1,\ldots$ and $a>0.$ They satisfy the discrete
orthogonality condition \cite{MR0372517}%
\[%
{\displaystyle\sum\limits_{j=0}^{\infty}}
\frac{a^{j}}{j!}C_{n}(j)C_{m}(j)=a^{-n}e^{a}n!\delta_{nm}.
\]
They are part of the Askey-scheme \cite{koekoek94askeyscheme} of
hypergeometric orthogonal polynomials:

\begin{center}
$%
\begin{array}
[c]{ccccccccc}%
_{4}F_{3} & \fbox{Wilson} & \  & \fbox{\ Racah} &  &  &  &  & \\
& \downarrow\quad\searrow &  & \downarrow\quad\searrow &  &  &  &  & \\
_{3}F_{2} & \ \fbox{$%
\begin{array}
[c]{c}%
\text{Continuous }\\
\text{dual Hahn}%
\end{array}
$} & \fbox{$%
\begin{array}
[c]{c}%
\text{Continuous }\\
\text{Hahn}%
\end{array}
$}\  & \fbox{Hahn} & \fbox{Dual Hahn} &  &  &  & \\
& \downarrow & \swarrow\quad\downarrow & \swarrow\quad\downarrow\quad\searrow
& \swarrow\quad\downarrow &  &  &  & \\
_{2}F_{1} & \fbox{$%
\begin{array}
[c]{c}%
\text{Meixner}\\
\text{Pollaczek}%
\end{array}
$} & \fbox{Jacobi} & \ \fbox{Meixner} & \fbox{Krawtchouk} &  &  &  & \\
& \searrow & \downarrow & \swarrow\quad\searrow & \downarrow &  &  &  & \\
_{1}F_{1} &  & \fbox{Laguerre} &  & \fbox{Charlier\ } & _{2}F_{0} &  &  & \\
&  & \searrow &  & \swarrow &  &  &  & \\
_{2}F_{0} &  &  & \fbox{Hermite} &  &  &  &  &
\end{array}
$
\end{center}

where the arrows indicate limit relations between the polynomials.

The Charlier polynomials have applications in quantum mechanics
\cite{MR1383073}, \cite{MR1985710}, \cite{MR1912384}, \cite{MR1381407},
difference equations \cite{MR1327166}, \cite{MR1677708}, teletraffic theory
\cite{MR0376131}$,$ \cite{MR954719}, generating functions \cite{MR863674},
\cite{MR0432949}, \cite{MR0279355}, and probability theory \cite{MR885141},
\cite{MR1693293}, \cite{MR1978097}, \cite{MR1662708}. The $q$-analogue of the
Charlier polynomials were studied in \cite{MR1947557}, \cite{MR1309153},
\cite{MR1418759} and \cite{MR1335793}. The generalized Charlier polynomials
were analyzed in \cite{MR1737084}, \cite{MR0382752}, \cite{MR1192598},
\cite{MR735537} and \cite{MR2063533}.

Asymptotics for the $L^{p}$-norms and information entropies of Charlier
polynomials were derived in \cite{MR1920124}. Bounds for their zeros were
obtained in \cite{MR2006607}. Asymptotic representations were established in
\cite{MR2026026} in terms of Hermite polynomials and in \cite{MR2068115} in
terms of Gamma functions. Some asymptotic estimates were computed in
\cite{MR1341423} from a representation of $C_{n}(x)$ in terms of Bell
polynomials. An asymptotic formula when $x<0$ was derived in \cite{MR778685}
using probabilistic methods.

In \cite{MR1606887}, Goh studied the asymptotic behavior of $C_{n}(x)$ for
large $n$ using an approximation of the Plancharel-Rotach type. A uniform
asymptotic expansion was derived in \cite{MR1297273} using the saddle-point
method. Asymptotic expansions were obtained in \cite{MR1857605} from a second
order linear differential satisfied by $C_{n}(x)$ in which $a$ is the
independent variable and $x$ is a parameter. 

In this paper we shall take a different approach and investigate the
asymptotic behavior of $C_{n}(x)$ as $n\rightarrow\infty,$ by using the limit
relation between the Krawtchouk polynomials $K_{n}(x)$ defined by
\begin{equation}
K_{n}(x)=K_{n}(x,p,N)=\,_{2}F_{1}\left(  \left.
\begin{array}
[c]{c}%
-n,-x\\
-N
\end{array}
\right\vert \frac{1}{p}\right)  ,\quad n=0,1,\ldots,N,\ 0\leq x\leq N,\ 0\leq
p\leq1 \label{K}%
\end{equation}
and the Charlier polynomials, namely%
\begin{equation}
\underset{N\rightarrow\infty}{\lim}K_{n}\left(  x,\frac{a}{N},N\right)
=C_{n}(x). \label{limit}%
\end{equation}
We shall use the asymptotic expansions derived in \cite{kraw} for the scaled
Krawtchouk polynomials $k_{n}(x)$, with%
\begin{equation}
k_{n}(x)=k_{n}(x,p,N)=\left(  -p\right)  ^{n}\binom{N}{n}K_{n}(x,p,N).
\label{kK}%
\end{equation}

A similar idea has been used in \cite{MR1985821} and \cite{MR1858318} to
obtain asymptotic approximations of several orthogonal polynomials of the
Askey-scheme in terms of Hermite and Laguerre polynomials.

\section{Preliminaries}

The following is the main result derived in \cite{kraw}.

\begin{theorem}
As $N\rightarrow\infty,$ $k_{n}(x,p,N)$ admits the following asymptotic
approximations (see Figure\ref{regions1}).

\begin{figure}[ptb]
\begin{center}
\rotatebox{270} {\resizebox{12cm}{!}{\includegraphics{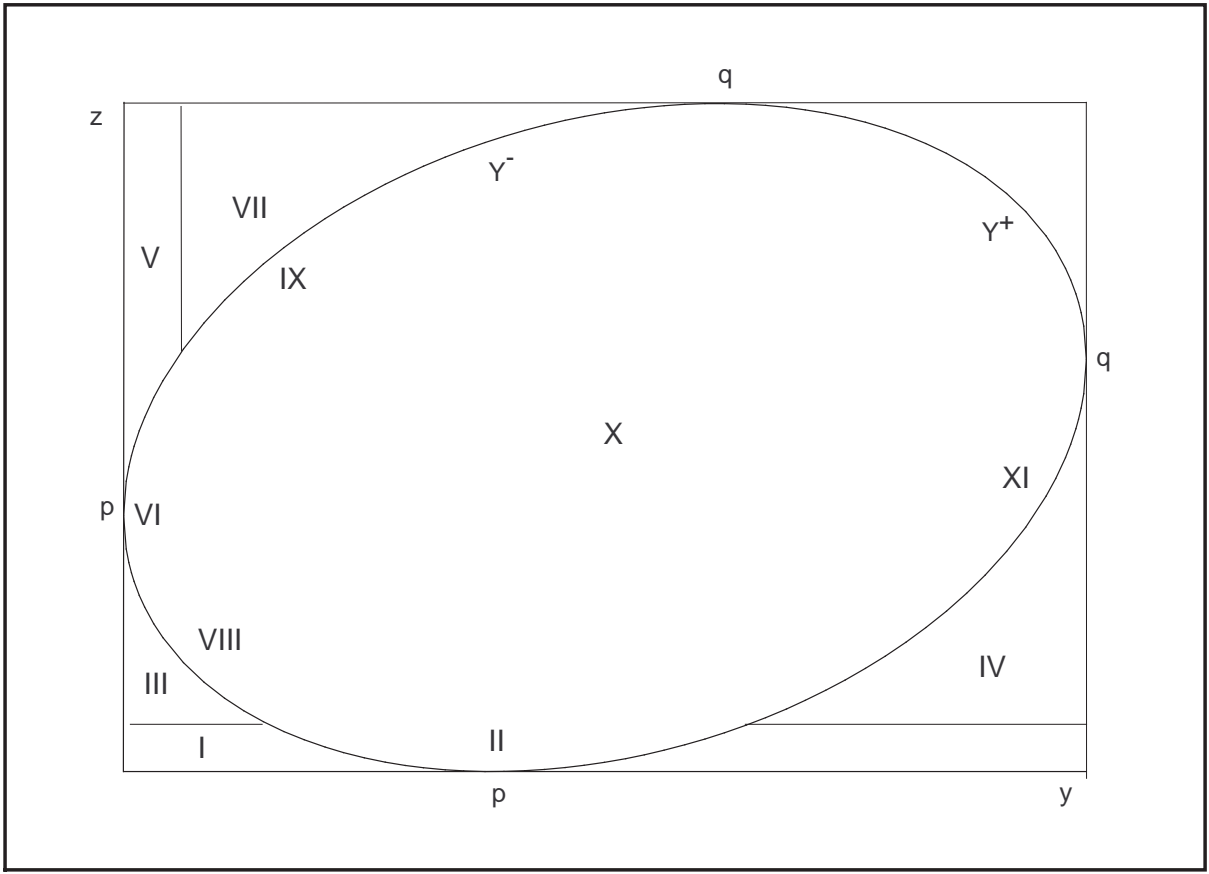}}}
\end{center}
\caption{A sketch of the different asymptotic regions for $k_{n}(x)$.}%
\label{regions1}%
\end{figure}

\begin{enumerate}
\item $n=O(1),$ $0\leq y\leq1,$ $y\not \approx p.$%
\begin{equation}
k_{n}(x)\sim k_{n}^{(1)}(y)=\frac{\varepsilon^{-n}}{n!}\left(  y-p\right)
^{n}\label{k1}%
\end{equation}
where
\[
\varepsilon=N^{-1},\quad x=\frac{y}{\varepsilon},\quad n=\frac{z}{\varepsilon
}\quad0\leq y,z\leq1.
\]

\item $n=O(1),$ $y\approx p,$ $y=p+\eta\sqrt{2pq\varepsilon},$ $\eta=O(1).$
\begin{equation}
k_{n}(x)\sim k_{n}^{\left(  2\right)  }(\eta)=\frac{\varepsilon^{-\frac{n}{2}%
}}{n!}\left(  \frac{pq}{2}\right)  ^{\frac{n}{2}}\mathrm{H}_{n}\left(
\eta\right)  ,\label{k2}%
\end{equation}
where $q=1-p$ and $\mathrm{H}_{n}\left(  \eta\right)  $ is the Hermite polynomial.

\item $0\leq y<Y^{-}(z),$ $0<z<p,$ where
\begin{equation}
Y^{\pm}(z)=p+\left(  q-p\right)  z\pm2zU_{0},\quad U_{0}(z)=\sqrt
{\frac{pq(1-z)}{z}}.\label{Ypm}%
\end{equation}%
\begin{equation}
k_{n}(x)\sim k^{(3)}(y,z)=\frac{\sqrt{\varepsilon}}{\sqrt{2\pi}}\exp\left[
\psi(y,z,U^{-})\varepsilon^{-1}\right]  L(z,U^{-}),\label{k3}%
\end{equation}
with%
\begin{equation}
\psi(y,z)=\left(  z-1\right)  \ln(U)+(1-y)\ln(U-p)+y\ln(U+q),\label{psi}%
\end{equation}%
\begin{equation}
L(y,z)=\sqrt{\frac{(U-p)(U+q)}{z\left[  U^{2}-\left(  U_{0}\right)
^{2}\right]  }}\label{L}%
\end{equation}
and%
\begin{equation}
U^{\pm}(y,z)=-\frac{1}{2}\left(  \frac{p-y}{z}+q-p\right)  \pm\frac{1}{2}%
\sqrt{\left(  \frac{p-y}{z}+q-p\right)  ^{2}-4\left(  U_{0}\right)  ^{2}%
}.\label{U}%
\end{equation}

\item $Y^{+}(z)<y\leq1,$ $0<z<q.$%
\begin{equation}
k_{n}(x)\sim k^{\left(  4\right)  }(y,z)=\frac{\sqrt{\varepsilon}}{\sqrt{2\pi
}}\exp\left[  \psi(y,z,U^{+})\varepsilon^{-1}\right]  L(z,U^{+}).\label{k4}%
\end{equation}

\item $x=O(1),$ $p<z<1.$%
\begin{gather}
k_{n}(x)\sim k^{(5)}\left(  x,z\right)  =\frac{\sqrt{\varepsilon}}{\sqrt{2\pi
}\sqrt{z\left(  1-z\right)  }}\cos(\pi x)\left(  \frac{z-p}{p}\right)
^{x}\exp\left[  \phi_{0}(z)\varepsilon^{-1}\right]  \label{k5}\\
-\frac{\varepsilon}{\pi}\frac{x}{z-p}\Gamma(x)\sin(\pi x)\left(
\frac{q\varepsilon}{z-p}\right)  ^{x}\exp\left[  \frac{(z-1)\ln(q)+\pi
\mathrm{i}z}{\varepsilon}\right]  \nonumber
\end{gather}
where%
\[
\phi_{0}(z)=(z-1)\ln(1-z)-z\ln(z)+z\ln(-p)
\]
and $\Gamma(x)$ is the Gamma function.

\item $x=O(1),$ $z\approx p,$ $z=p-u\sqrt{pq\varepsilon},$ $u=O(1).$%
\begin{align}
k_{n}(x) &  \sim k^{(6)}(x,u)=\frac{\sqrt{\varepsilon}}{\sqrt{2\pi pq}}\left[
\sqrt{\frac{q\varepsilon}{p}}\right]  ^{x}\mathrm{D}_{x}(u)\label{k6}\\
&  \times\exp\left[  \frac{\pi\mathrm{i}p-q\ln\left(  q\right)  }{\varepsilon
}+\frac{u\sqrt{pq}\pi\mathrm{i-}u\sqrt{pq}\ln\left(  q\right)  }%
{\sqrt{\varepsilon}}-\frac{u^{2}}{4}\right]  ,\nonumber
\end{align}
where $\mathrm{D}_{x}(u)$ is the parabolic cylinder function.

\item $0\ll y<Y^{-}(z),$ $p<z<1.$%
\begin{equation}
k_{n}(x)\sim k^{\left(  7\right)  }(y,z)=\exp\left(  \frac{\pi\mathrm{i}%
y}{\varepsilon}\right)  \left[  \cos\left(  \frac{\pi y}{\varepsilon}\right)
k^{(4)}(y,z)+2\mathrm{i}\sin\left(  \frac{\pi y}{\varepsilon}\right)
k^{(3)}(y,z)\right]  \label{k7}%
\end{equation}

\item $y\approx Y^{-}(z),$ $0<z<p,$ $y=Y^{-}(z)-\beta\varepsilon^{2/3},$
$\beta=O(1)$.%
\begin{equation}
k_{n}(x)\sim k^{(8)}(\beta,z)=\varepsilon^{\frac{1}{3}}\exp\left[  \psi
_{0}(z)\varepsilon^{-1}+\ln\left(  \frac{U_{0}+p}{U_{0}-q}\right)
\beta\varepsilon^{-\frac{1}{3}}\right]  \mathrm{Ai}\left(  \Theta^{^{\frac
{2}{3}}}\beta\right)  \frac{\Theta^{-\frac{1}{3}}}{\sqrt{zU_{0}}},\label{k8}%
\end{equation}
where%
\begin{equation}
\psi_{0}(z)=z\pi\mathrm{i}+(z-1)\ln\left(  U_{0}\right)  +Y^{-}(z)\ln\left(
U_{0}-q\right)  +\left[  1-Y^{-}(z)\right]  \ln\left(  U_{0}+p\right)
,\label{psi0}%
\end{equation}%
\begin{equation}
\Theta(z)=\sqrt{\frac{U_{0}}{z}}\frac{1}{\left(  U_{0}+p\right)  \left(
U_{0}-q\right)  }\label{Theta}%
\end{equation}
and $\mathrm{Ai}\left(  \cdot\right)  $ is the Airy function.

\item $y\approx Y^{-}(z),$ $p<z<1,\quad y=Y^{-}(z)-\beta\varepsilon^{\frac
{2}{3}},\ \beta=O(1).$%
\begin{align}
k_{n}(x) &  \sim k^{(9)}(\beta,z)=\varepsilon^{\frac{1}{3}}\exp\left[
\psi_{0}(z)\varepsilon^{-1}+\ln\left(  \frac{U_{0}+p}{U_{0}-q}\right)
\beta\varepsilon^{-\frac{1}{3}}\right]  \label{k9}\\
&  \times\frac{1}{2}\frac{\vartheta^{-\frac{1}{3}}}{\sqrt{zU_{0}}}\left[
\lambda^{+}(\beta,z)\mathrm{Ai}\left(  \vartheta^{^{\frac{2}{3}}}\beta\right)
+\mathrm{i}\lambda^{-}(\beta,z)\mathrm{Bi}\left(  \vartheta^{^{\frac{2}{3}}%
}\beta\right)  \right]  ,\nonumber
\end{align}
where $\vartheta(z)=-\Theta(z),$%
\begin{equation}
\lambda^{\pm}(\beta,z)=\exp\left\{  2\pi\mathrm{i}\left[  Y^{-}\left(
z\right)  -\beta\varepsilon^{\frac{2}{3}}\right]  \varepsilon^{-1}\right\}
\pm1.\label{lambda}%
\end{equation}
and $\mathrm{Ai}\left(  \cdot\right)  ,\mathrm{Bi}\left(  \cdot\right)  $ are
the Airy functions.

\item $Y^{-}(z)<y<Y^{+}(z),$ $0<z<1.$%
\begin{equation}
k_{n}(x)\sim k^{(10)}(y,z)=k^{(3)}(y,z)+k^{(4)}(y,z).\label{k10}%
\end{equation}

\item $y\approx Y^{+}(z),$ $0<z<q,\quad y=Y^{+}(z)+\alpha\varepsilon^{\frac
{2}{3}},\quad\alpha=O(1).$%
\begin{equation}
k_{n}(x)\sim k^{(11)}(\alpha,z)=\varepsilon^{\frac{1}{3}}\exp\left[  \psi
_{1}(z)\varepsilon^{-1}+\ln\left(  \frac{U_{0}+q}{U_{0}-p}\right)
\alpha\varepsilon^{-\frac{1}{3}}\right]  \mathrm{Ai}\left[  \left(  \Theta
_{1}\right)  ^{\frac{2}{3}}\alpha\right]  \frac{\left(  \Theta_{1}\right)
^{-\frac{1}{3}}}{\sqrt{zU_{0}}},\label{k11}%
\end{equation}
where%
\begin{equation}
\psi_{1}(z)=(z-1)\ln\left(  U_{0}\right)  +Y^{+}(z)\ln\left(  U_{0}+q\right)
+\left[  1-Y^{+}(z)\right]  \ln\left(  U_{0}-p\right)  ,\label{psi1}%
\end{equation}%
\begin{equation}
\Theta_{1}(z)=\sqrt{\frac{U_{0}}{z}}\frac{1}{\left(  U_{0}-p\right)  \left(
U_{0}+q\right)  }\label{Theta1}%
\end{equation}

\end{enumerate}
\end{theorem}

In order to obtain the corresponding asymptotic expansions for $K_{n}(x),$ we
need to derive asymptotic formulas for $\left(  -p\right)  ^{n}\binom{N}{n}$
in the different regions of Theorem 1.

The following lemma follows immediately from Stirling's formula
\cite{MR94b:00012}%
\begin{equation}
\Gamma(x)\sim\sqrt{\frac{2\pi}{x}}x^{x}e^{-x},\quad x\rightarrow\infty.
\label{stirling}%
\end{equation}

\begin{lemma}
As $N\rightarrow\infty,$ we have the following asymptotic approximations:

\begin{enumerate}
\item
\begin{equation}
\left(  -p\right)  ^{n}\binom{N}{n}\sim\frac{1}{n!}\left(  -p\right)
^{n}\varepsilon^{-n},\quad n=O(1). \label{Bi1}%
\end{equation}

\item
\begin{equation}
\left(  -p\right)  ^{n}\binom{N}{n}\sim\frac{\sqrt{\varepsilon}}{\sqrt{2\pi
}\sqrt{z\left(  1-z\right)  }}\exp\left[  \phi(z)\varepsilon^{-1}\right]
,\quad n=z\varepsilon^{-1} \label{Bi2}%
\end{equation}
where%
\begin{equation}
\phi(z)=z\ln(p)+z\pi\mathrm{i}+(z-1)\ln(1-z)-z\ln(z). \label{phi}%
\end{equation}

\item
\begin{equation}
\left(  -p\right)  ^{n}\binom{N}{n}\sim\frac{\sqrt{\varepsilon}}{\sqrt{2\pi
pq}}\exp\left[  \phi_{1}(u)\right]  ,\quad n=p\varepsilon^{-1}-u\sqrt
{\frac{pq}{\varepsilon}},\quad u=O(1) \label{Bi3}%
\end{equation}
where%
\begin{equation}
\phi_{1}(u)=\left[  \pi\mathrm{i}p-q\ln(q)\right]  \varepsilon^{-1}-u\sqrt
{pq}\left[  \ln(q)+\pi\mathrm{i}\right]  \varepsilon^{-\frac{1}{2}}-\frac
{1}{2}u^{2}. \label{phi1}%
\end{equation}

\end{enumerate}
\end{lemma}

\section{Limit analysis}

\begin{figure}[ptb]
\begin{center}
\rotatebox{270} {\resizebox{12cm}{!}{\includegraphics{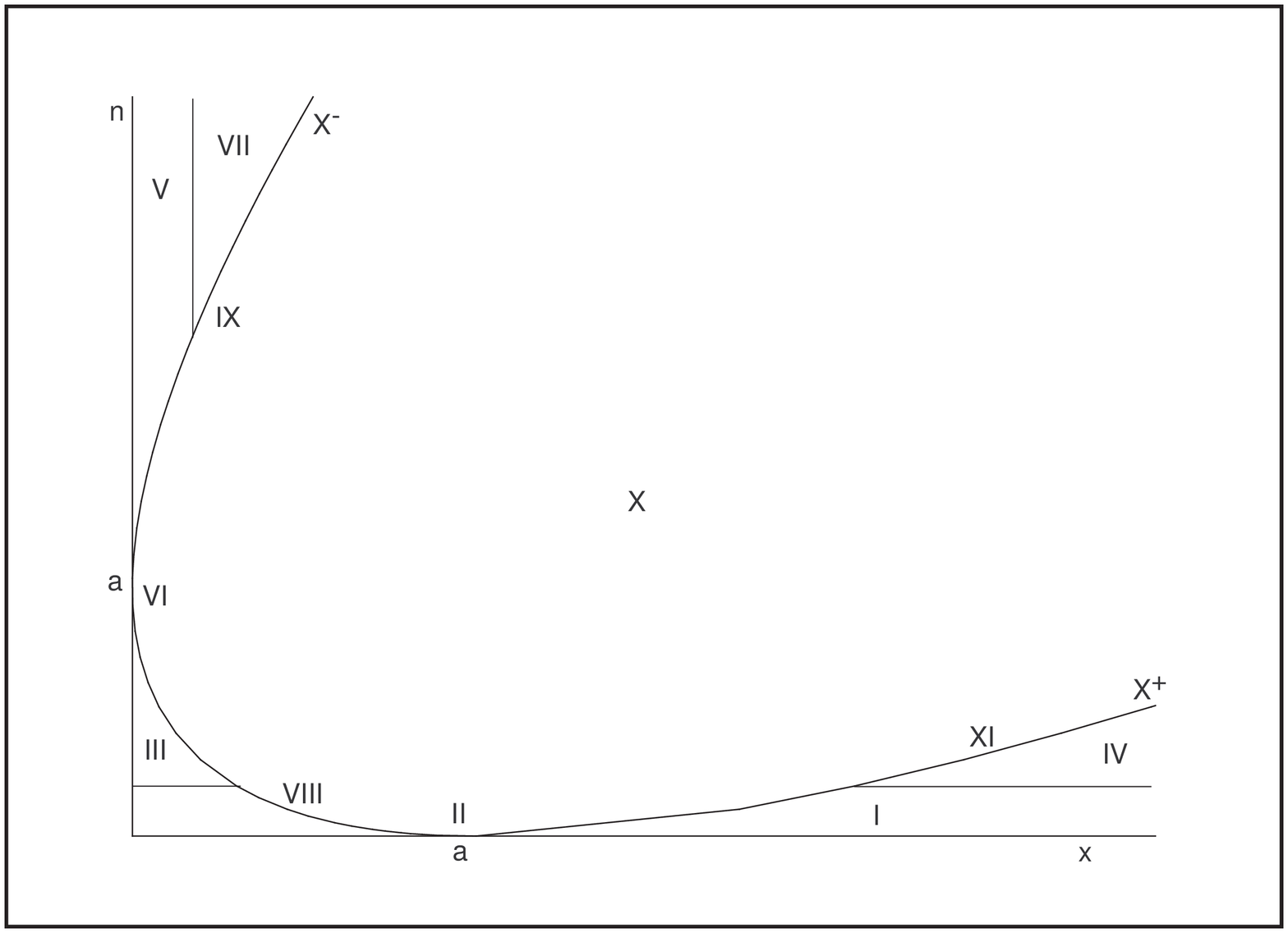}}}
\end{center}
\caption{A sketch of the different asymptotic regions for $C_{n}(x)$.}%
\label{regions2}%
\end{figure}

Setting
\begin{equation}
y=x\varepsilon,\quad z=n\varepsilon,\quad q=1-p,\quad p=a\varepsilon
\label{yzq}%
\end{equation}
in (\ref{Ypm}) and letting $\varepsilon\rightarrow0,$ we obtain%
\begin{equation}
Y^{\pm}(z)\rightarrow X^{\pm}(a,n)=\left(  \sqrt{n}\pm\sqrt{a}\right)
^{2},\quad U_{0}(z)\rightarrow\sqrt{\frac{a}{n}}.\label{Xpm}%
\end{equation}
Hence, the eleven regions of Theorem 1 transform into the following regions
(see Figure \ref{regions2}).

\begin{enumerate}
\item Region I

From (\ref{k1}) and (\ref{Bi1}) we have for $n=O(1)$
\begin{equation}
K_{n}(x)\sim\left(  1-\frac{y}{p}\right)  ^{n}.\label{K1}%
\end{equation}
Thus,%
\begin{equation}
C_{n}(x)\simeq\left(  1-\frac{x}{a}\right)  ^{n}.\label{G1}%
\end{equation}
The formula above is exact for $n=0,1$ and in the limit as $a\rightarrow
\infty$ we have%
\[
\underset{a\rightarrow\infty}{\lim}C_{n}(av)=\left(  1-v\right)  ^{n}%
=~_{1}F_{0}\left(  \left.
\begin{array}
[c]{c}%
-n\\
-
\end{array}
\right\vert \ v\right)  .
\]

\item Region II

From (\ref{k2}) and (\ref{Bi1}) we have for $x=p\varepsilon^{-1}+\eta
\sqrt{\frac{2pq}{\varepsilon}},\quad\eta=O(1)$
\begin{equation}
K_{n}(x)\sim\left(  -1\right)  ^{n}\left(  \frac{\varepsilon q}{2p}\right)
^{\frac{n}{2}}\mathrm{H}_{n}\left(  \eta\right)  . \label{K2}%
\end{equation}
Hence,%
\begin{equation}
C_{n}(x)\simeq\left(  -1\right)  ^{n}\left(  2a\right)  ^{-\frac{n}{2}%
}\mathrm{H}_{n}\left(  \eta\right)  ,\quad x=a+\eta\sqrt{2a}. \label{G2}%
\end{equation}

Equation (\ref{G2}) is exact for $n=0,1$ and in the limit as $a\rightarrow
\infty$ we have%
\[
\underset{a\rightarrow\infty}{\lim}\left(  -1\right)  ^{n}\left(  2a\right)
^{\frac{n}{2}}C_{n}(a+\eta\sqrt{2a})=\mathrm{H}_{n}\left(  \eta\right)
\]
which is equation 2.12.1 in.

\item Region III

From (\ref{k3}), (\ref{Bi2}) and (\ref{phi}) we have for $0\leq y<Y^{-}%
(z),\quad0<z<p$
\begin{equation}
K_{n}(x)\sim K^{\left(  3\right)  }(y,z)=\exp\left[  \frac{\psi\left(
y,z,U^{-}\right)  -\phi(z)}{\varepsilon}\right]  G\left(  z,U^{-}\right)
,\label{K3}%
\end{equation}
where%
\begin{equation}
G(z,U)=\sqrt{\frac{(1-z)(U-p)(U+q)}{U^{2}-Uo^{2}}}.\label{G}%
\end{equation}
From (\ref{psi}), (\ref{U}) and (\ref{G}) we obtain%
\[
\frac{\psi\left(  y,z,U^{-}\right)  -\phi(z)}{\varepsilon}\rightarrow\Psi
_{3}(x),
\]
where%
\begin{equation}
\Psi_{3}(x)=x\ln\left(  \frac{a+x-n+\Delta}{2a}\right)  +n\ln\left(
\frac{a-x+n+\Delta}{2a}\right)  +\frac{1}{2}\left(  a-x-n+\Delta\right)
\label{Psi3}%
\end{equation}
and%
\begin{equation}
G\left(  z,U^{-}\right)  \rightarrow L_{3}(x)\equiv\sqrt{\frac{a-x-n+\Delta
}{2\Delta}},\label{L3}%
\end{equation}
for $0\leq x<X^{-},\quad0<n<a,$ with%
\begin{equation}
\Delta(a,n,x)=\sqrt{a^{2}-2a(x+n)+(x-n)^{2}}.\label{Delta}%
\end{equation}
Thus,%
\begin{equation}
C_{n}(x)\sim F_{3}(x)=\exp\left[  \Psi_{3}(x)\right]  L_{3}(x),\quad0\leq
x<X^{-},\quad0<n<a.\label{C3}%
\end{equation}
Figure \ref{RegionIII} shows the accuracy of the approximation (\ref{C3}) with
$n=30$ and $a=50.165184$ in the range $-3<x<X^{-}$.

\begin{figure}[ptb]
\begin{center}
\rotatebox{270} {\resizebox{12cm}{!}{\includegraphics{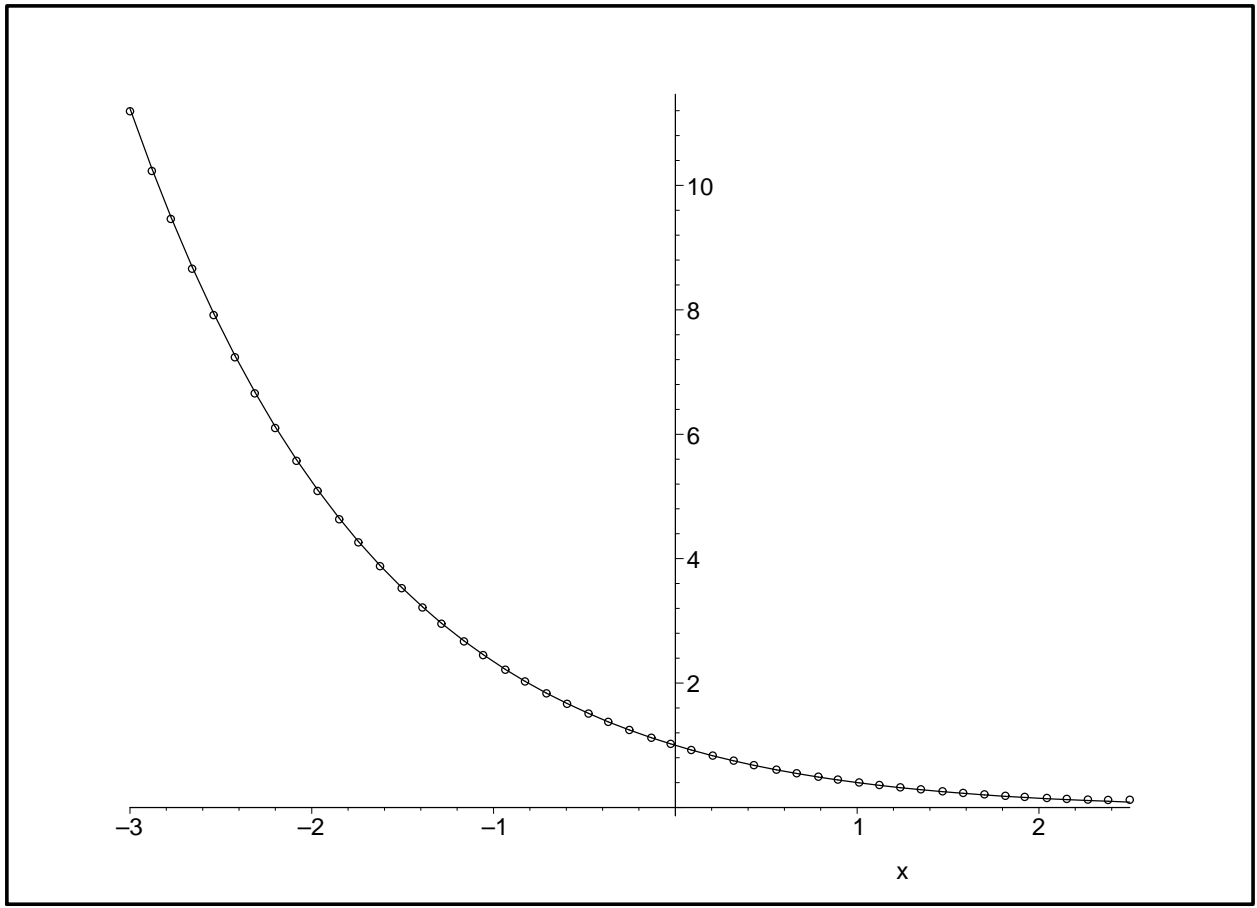}}}
\end{center}
\caption{A comparison of $C_{n}(x)$ (solid curve) and $F_{3}(x)$ (ooo) for $n=30$ with $a=50.165184$.}%
\label{RegionIII}%
\end{figure}

\begin{remark}
Although we have not shown proof of it, the approximation (\ref{C3}) is in
fact also valid for $x<0$ and $n\geq0$ (see Figure \ref{RegionVneg}).

\begin{figure}[ptb]
\begin{center}
\rotatebox{270} {\resizebox{12cm}{!}{\includegraphics{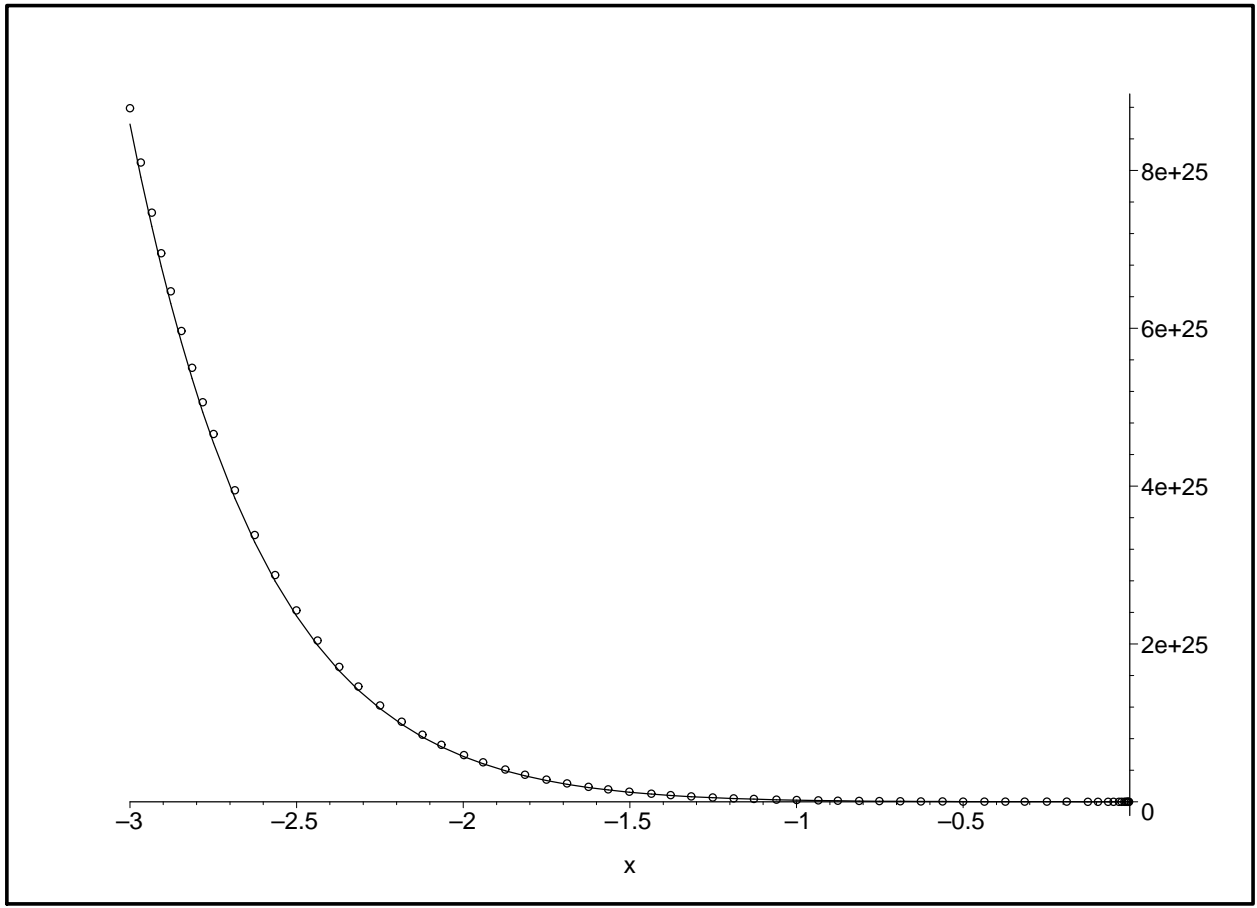}}}
\end{center}
\caption{A comparison of $C_{n}(x)$ (solid curve) and $F_{3}(x)$ (ooo) for $n=30$ with $a=2.165184$.}%
\label{RegionVneg}%
\end{figure}
\end{remark}

\item Region IV

From (\ref{k4}), (\ref{Bi2}) and (\ref{phi}) we have for $Y^{+}(z)<y\leq
1,\quad0<z<q$
\begin{equation}
K_{n}(x)\sim K^{\left(  4\right)  }(y,z)=\exp\left[  \frac{\psi\left(
y,z,U^{+}\right)  -\phi(z)}{\varepsilon}\right]  G\left(  z,U^{+}\right)
.\label{K4}%
\end{equation}
From (\ref{psi}), (\ref{U}) and (\ref{G}) we obtain%
\[
\frac{\psi\left(  y,z,U^{+}\right)  -\phi(z)}{\varepsilon}\rightarrow\Psi
_{4}(x)-n\pi\mathrm{i},
\]
where%
\begin{equation}
\Psi_{4}(x)=x\ln\left(  \frac{a+x-n-\Delta}{2a}\right)  +n\ln\left(
\frac{x-a-n+\Delta}{2a}\right)  +\frac{1}{2}\left(  a-x-n+\Delta\right)
\label{Psi4}%
\end{equation}
and%
\begin{equation}
G\left(  z,U^{-}\right)  \rightarrow L_{4}(x)\equiv\sqrt{\frac{x-a+n+\Delta
}{2\Delta}},\label{L4}%
\end{equation}
for $X^{+}<x$. Therefore,%
\begin{equation}
C_{n}(x)\sim F_{4}(x)=\left(  -1\right)  ^{n}\exp\left[  \Psi_{4}(x)\right]
L_{4}(x),\quad X^{+}<x.\label{C4}%
\end{equation}

Figure \ref{RegionIV} shows the accuracy of the approximation (\ref{C4}) with
$n=30$ and $a=2.165184$ in the range $X^{+}<x<\infty$.

\begin{figure}[ptb]
\begin{center}
\rotatebox{270} {\resizebox{12cm}{!}{\includegraphics{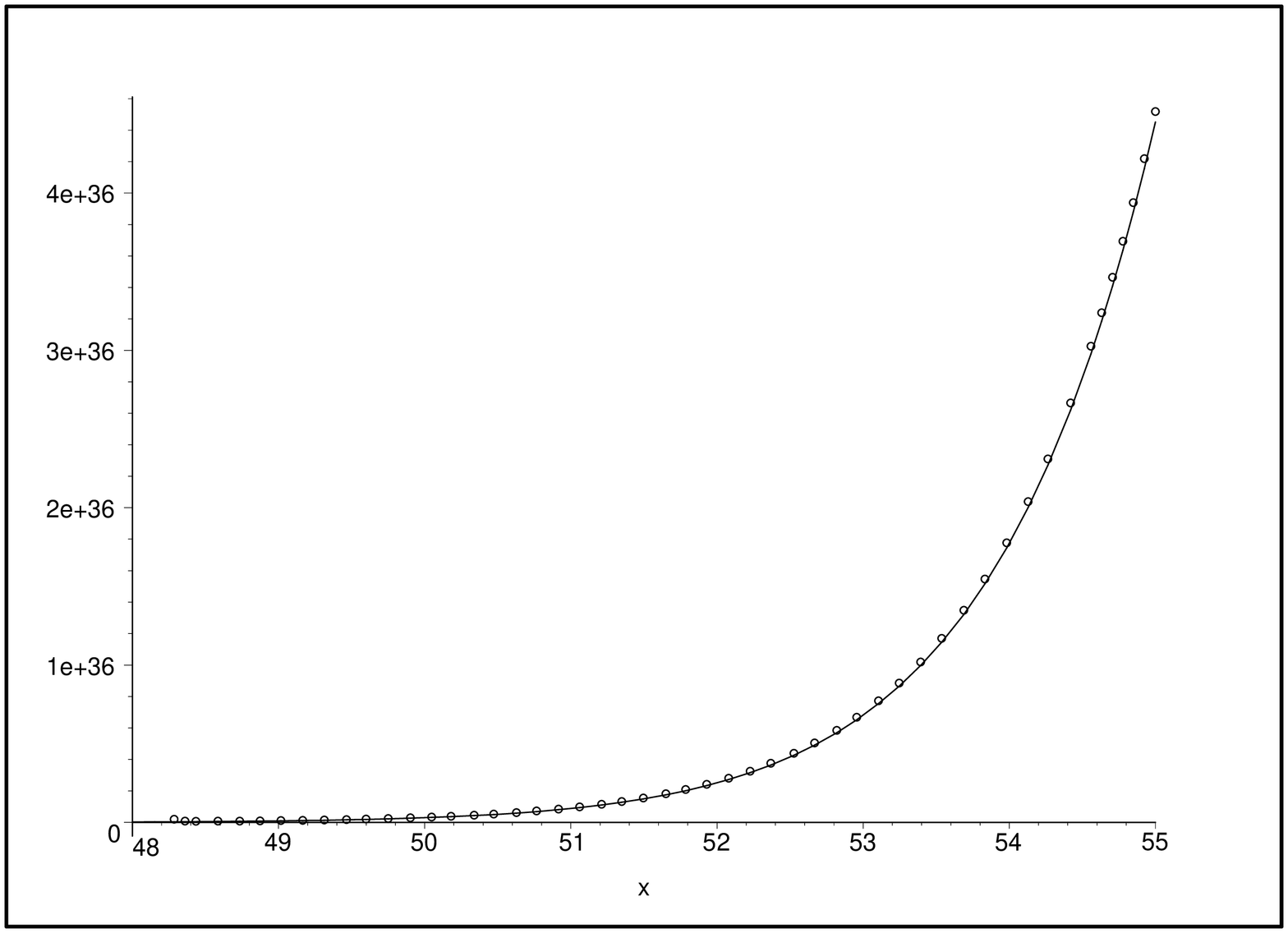}}}
\end{center}
\caption{A comparison of $C_{n}(x)$ (solid curve) and $F_{4}(x)$ (ooo) for $n=30$ with $a=2.165184$.}%
\label{RegionIV}%
\end{figure}

\item Region V

From (\ref{k5}), (\ref{Bi2}) and (\ref{phi}) we have for $x=O(1),\quad p<z<1$
\begin{gather}
K_{n}(x)\sim\exp\left[  x\ln\left(  \frac{z}{p}-1\right)  \right]  \cos\left(
\pi x\right)  -\frac{\sqrt{\varepsilon}}{z-p}\sqrt{\frac{2}{\pi}z\left(
1-z\right)  }\Gamma\left(  x+1\right)  \sin\left(  \pi x\right)  \label{K5}\\
\times\exp\left[  \frac{(1-z)\ln\left(  \frac{1-z}{q}\right)  +z\ln\left(
\frac{z}{p}\right)  }{\varepsilon}+x\ln\left(  \frac{\varepsilon q}%
{z-p}\right)  \right]  .\nonumber
\end{gather}
Hence,%
\begin{gather}
C_{n}(x)\sim F_{5}(x)=\exp\left[  x\ln\left(  \frac{n}{a}-1\right)  \right]
\cos\left(  \pi x\right)  -\sqrt{n}\sqrt{\frac{2}{\pi}}\Gamma\left(
x+1\right)  \sin\left(  \pi x\right)  \label{C5}\\
\times\exp\left[  n\ln\left(  \frac{n}{a}\right)  -(x+1)\ln(n-a)+a-n\right]
,\quad x\approx0,\quad n>a.\nonumber
\end{gather}

Figure \ref{RegionV} shows the accuracy of the approximation (\ref{C5}) with
$n=30$, $a=2.165184$ and $x\approx0$.

\begin{figure}[ptb]
\begin{center}
\rotatebox{270} {\resizebox{12cm}{!}{\includegraphics{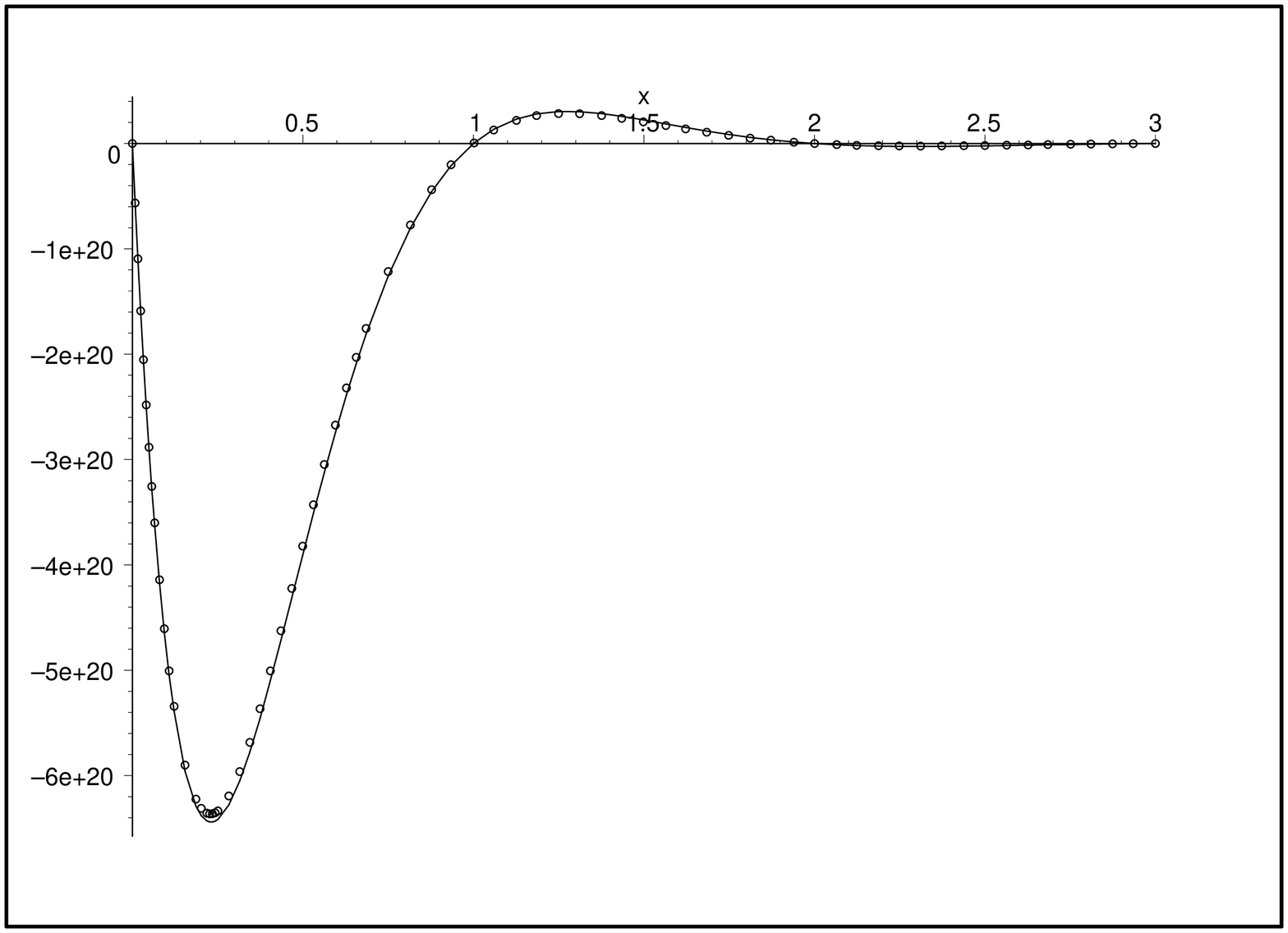}}}
\end{center}
\caption{A comparison of $C_{n}(x)$ (solid curve) and $F_{5}(x)$ (ooo) for $n=30$ with $a=2.165184$.}%
\label{RegionV}%
\end{figure}

\item Region VI

From (\ref{k6}) and (\ref{Bi3})-(\ref{phi1}) we have for $x=O(1),\quad
z=p-u\sqrt{pq\varepsilon},\quad u=O(1)$
\begin{equation}
K_{n}(x)\sim\exp\left[  \frac{x}{2}\ln\left(  \frac{q\varepsilon}{p}\right)
+\frac{u^{2}}{4}\right]  \mathrm{D}_{x}\left(  u\right)  .\label{K6}%
\end{equation}
Thus,%
\begin{equation}
C_{n}(x)\sim F_{6}(x)=\exp\left[  -\frac{x}{2}\ln\left(  a\right)
+\frac{u^{2}}{4}\right]  \mathrm{D}_{x}\left(  u\right)  ,\quad x\approx
0,\quad n=a-u\sqrt{a}.\label{C6}%
\end{equation}

Figure \ref{RegionVI} shows the accuracy of the approximation (\ref{C6}) with
$n=30$, $a=30.165184$ ($n\approx a$) and $x\approx0$.

\begin{figure}[ptb]
\begin{center}
\rotatebox{270} {\resizebox{12cm}{!}{\includegraphics{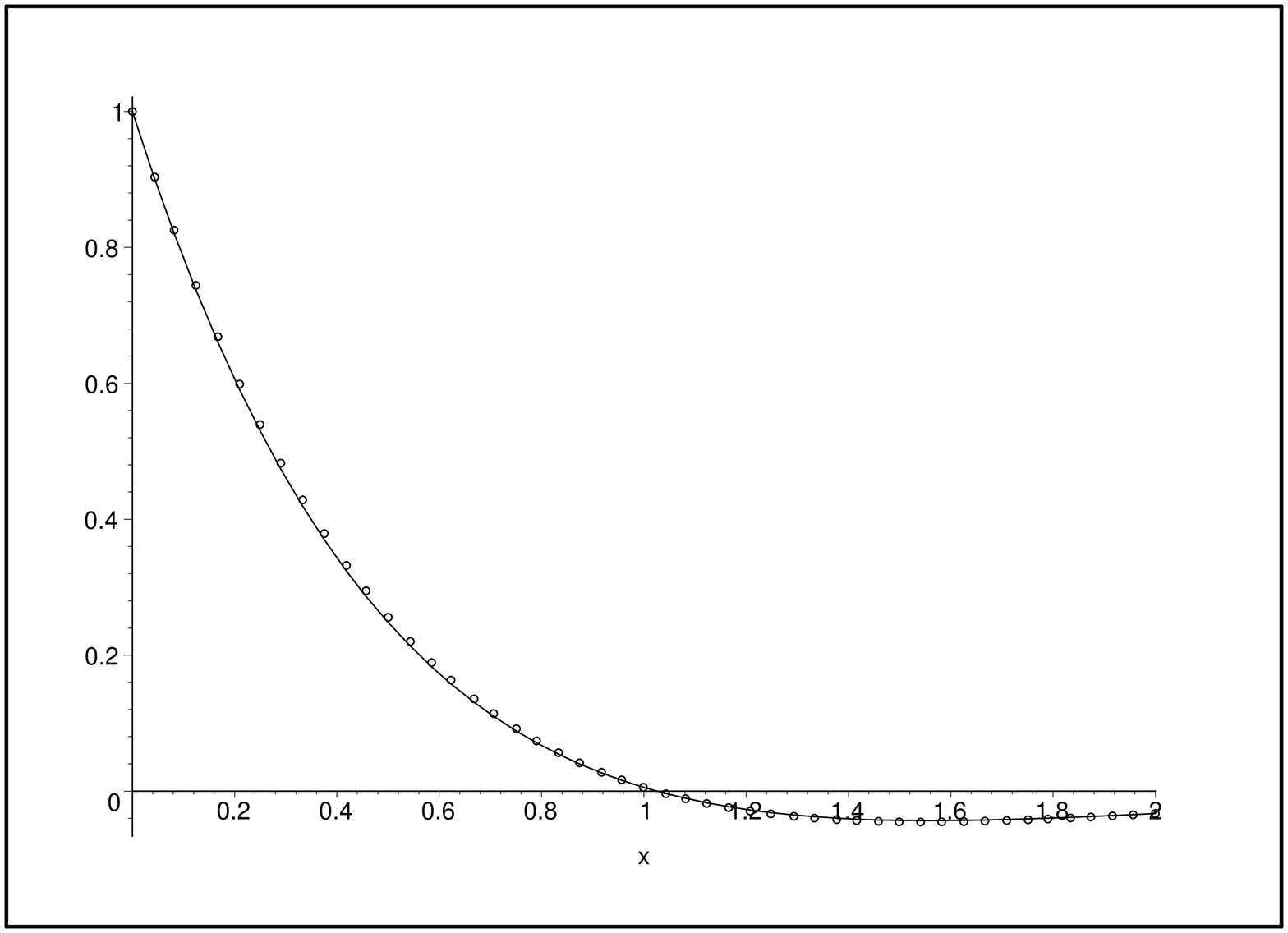}}}
\end{center}
\caption{A comparison of $C_{n}(x)$ (solid curve) and $F_{6}(x)$ (ooo) for $n=30$ with $a=30.165184$.}%
\label{RegionVI}%
\end{figure}

\item Region VII

From (\ref{k7}), (\ref{K3}) and (\ref{K4}) we have for $0\ll y<Y^{-}(z),$
$p<z<1$%
\begin{equation}
K_{n}(x)\sim\exp\left(  \frac{\pi\mathrm{i}y}{\varepsilon}\right)  \left[
\cos\left(  \frac{\pi y}{\varepsilon}\right)  K^{(4)}(y,z)+2\mathrm{i}%
\sin\left(  \frac{\pi y}{\varepsilon}\right)  K^{(3)}(y,z)\right]  .\label{K7}%
\end{equation}
Therefore%
\[
C_{n}(x)\sim\exp\left(  \pi\mathrm{i}x\right)  \left[  \cos\left(  \pi
x\right)  C_{n}^{(4)}(x)+2\mathrm{i}\sin\left(  \pi x\right)  C_{n}%
^{(3)}(x)\right]
\]
for $0\ll x<X^{-},\ n>a,$ which we can rewrite as%
\begin{align}
C_{n}(x) &  \sim F_{7}(x)=\exp\left[  x\ln\left(  \frac{n-a-x+\Delta}%
{2a}\right)  +n\ln\left(  \frac{a+n-x-\Delta}{2a}\right)  +\frac{1}{2}\left(
a-x-n+\Delta\right)  \right]  \nonumber\\
&  \times\cos\left(  \pi x\right)  \sqrt{\frac{x+n-a+\Delta}{2\Delta}}%
-2\sin\left(  \pi x\right)  \sqrt{\frac{x+n-a-\Delta}{2\Delta}}\label{C7}\\
&  \times\exp\left[  x\ln\left(  \frac{n-a-x-\Delta}{2a}\right)  +n\ln\left(
\frac{a+n-x+\Delta}{2a}\right)  +\frac{1}{2}\left(  a-x-n+\Delta\right)
\right]  .\nonumber
\end{align}

Figure \ref{RegionVII} shows the accuracy of the approximation (\ref{C7}) with
$n=30$ and $a=2.165184$ in the range $0\ll x<X^{-}$.

\begin{figure}[ptb]
\begin{center}
\rotatebox{270} {\resizebox{12cm}{!}{\includegraphics{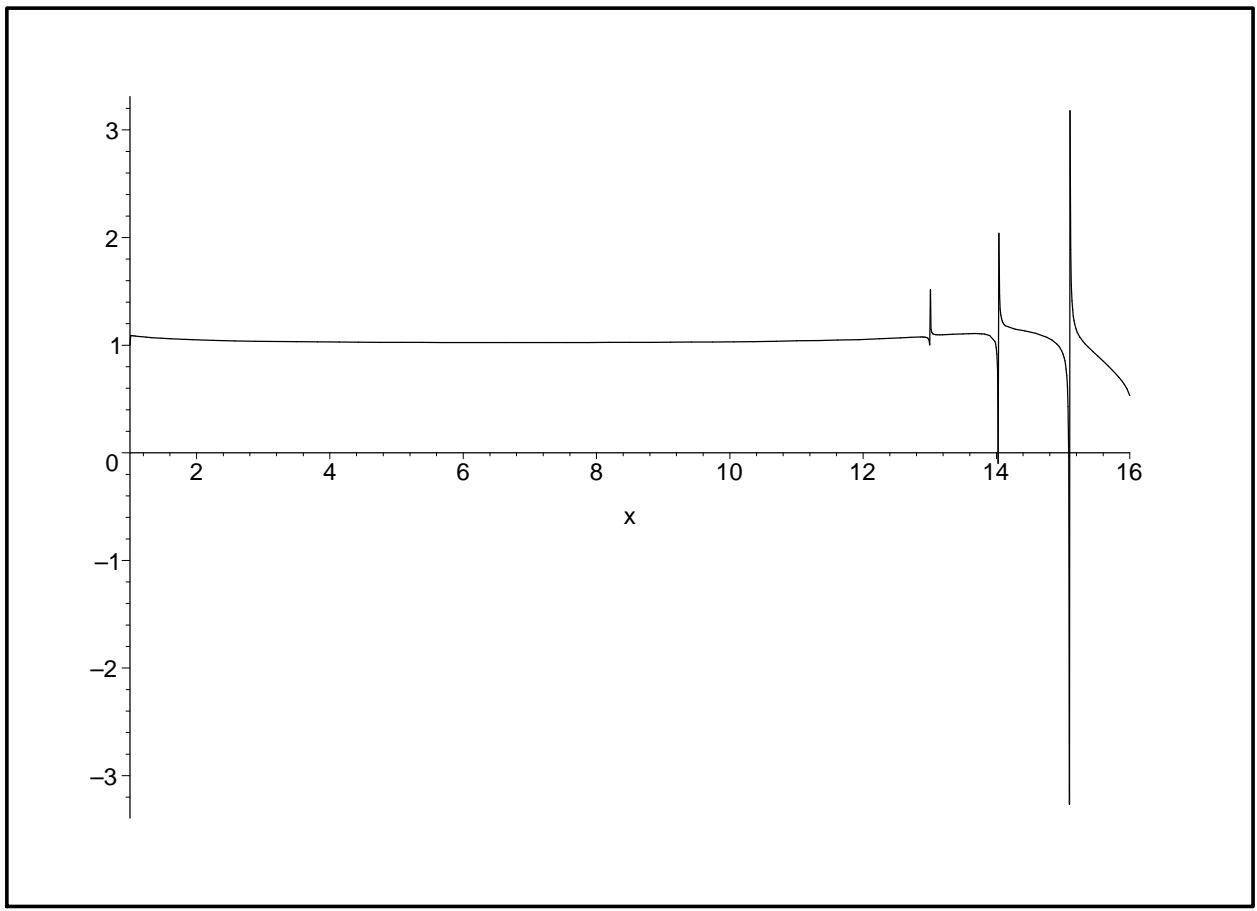}}}
\end{center}
\caption{A sketch of $C_{n}(x)/F_{7}(x)$ for $n=30$ with $a=2.165184$. The vertical lines are due to the discontinuities at $x\simeq 13,14$ and $15$.}%
\label{RegionVII}%
\end{figure}

\item Region VIII

From (\ref{k8}), (\ref{Bi2}) and (\ref{phi}) we have for \ $y\approx
Y^{-}(z),$ $0<z<p,$ $y=Y^{-}(z)-\beta\varepsilon^{2/3},$ $\beta=O(1)$%
\begin{align}
K_{n}(x) &  \sim\varepsilon^{-\frac{1}{6}}\exp\left[  \frac{\psi_{0}%
(z)-\phi(z)}{\varepsilon}+\ln\left(  \frac{U_{0}+p}{U_{0}-q}\right)
\beta\varepsilon^{-\frac{1}{3}}\right]  \label{K8}\\
&  \times\sqrt{2\pi}\left[  \frac{z\left(  1-z\right)  }{pq}\right]
^{\frac{1}{4}}\mathrm{Ai}\left(  \Theta^{^{\frac{2}{3}}}\beta\right)
\Theta^{-\frac{1}{3}}.\nonumber
\end{align}
Since
\begin{equation}
\beta=\left[  Y^{-}(z)-y\right]  \varepsilon^{-\frac{2}{3}}\label{beta}%
\end{equation}
we have from (\ref{Xpm})%
\begin{equation}
\beta\varepsilon^{-\frac{1}{3}}\rightarrow X^{-}-x.\label{betalim}%
\end{equation}
From (\ref{psi0}), (\ref{phi}) and (\ref{betalim}) we get%
\[
\frac{\psi_{0}(z)-\phi(z)}{\varepsilon}+\ln\left(  \frac{U_{0}+p}{U_{0}%
-q}\right)  \beta\varepsilon^{-\frac{1}{3}}\rightarrow\frac{1}{2}n\ln\left(
\frac{n}{a}\right)  +x\ln\left(  1-\sqrt{\frac{n}{a}}\right)  +\sqrt{an}%
-\sqrt{n}.
\]
From (\ref{Theta}) and (\ref{betalim}) we obtain%
\begin{equation}
\varepsilon^{-\frac{1}{6}}\sqrt{2\pi}\left[  \frac{z\left(  1-z\right)  }%
{pq}\right]  ^{\frac{1}{4}}\Theta^{-\frac{1}{3}}\rightarrow\sqrt{2\pi}\left(
\frac{n}{a}\right)  ^{\frac{1}{6}}\left(  \sqrt{a}-\sqrt{n}\right)  ^{\frac
{1}{3}}\label{limfac}%
\end{equation}
and%
\begin{equation}
\Theta^{^{\frac{2}{3}}}\beta\rightarrow\left(  \frac{n}{a}\right)  ^{\frac
{1}{6}}\frac{\left(  X^{-}-x\right)  }{\left(  \sqrt{a}-\sqrt{n}\right)
^{\frac{2}{3}}}.\label{limbeta}%
\end{equation}
Therefore,%
\begin{align}
C_{n}(x) &  \sim F_{8}(x)=\sqrt{2\pi}\left(  \frac{n}{a}\right)  ^{\frac{1}%
{6}}\left(  \sqrt{a}-\sqrt{n}\right)  ^{\frac{1}{3}}\mathrm{Ai}\left[  \left(
\frac{n}{a}\right)  ^{\frac{1}{6}}\frac{\left(  X^{-}-x\right)  }{\left(
\sqrt{a}-\sqrt{n}\right)  ^{\frac{2}{3}}}\right]  \label{C8}\\
&  \times\exp\left[  \frac{1}{2}n\ln\left(  \frac{n}{a}\right)  +x\ln\left(
1-\sqrt{\frac{n}{a}}\right)  +\sqrt{an}-\sqrt{n}\right]  \nonumber
\end{align}
for $x\approx X^{-},\ 0<n<a.$

\item Region IX

From (\ref{k9}), (\ref{lambda}), (\ref{Bi2}) and (\ref{phi}) we have for
$y\approx Y^{-}(z),$ $p<z<1,\quad y=Y^{-}(z)-\beta\varepsilon^{\frac{2}{3}%
},\ \beta=O(1)$%
\begin{align*}
K_{n}(x)  &  \sim\varepsilon^{-\frac{1}{6}}\exp\left[  \frac{\psi_{0}%
(z)-\phi(z)}{\varepsilon}+\ln\left(  \frac{U_{0}+p}{U_{0}-q}\right)
\beta\varepsilon^{-\frac{1}{3}}\right] \\
&  \times\sqrt{2\pi}\left[  \frac{z\left(  1-z\right)  }{pq}\right]
^{\frac{1}{4}}\vartheta^{-\frac{1}{3}}\left[  \lambda^{+}(\beta,z)\mathrm{Ai}%
\left(  \vartheta^{^{\frac{2}{3}}}\beta\right)  +\mathrm{i}\lambda^{-}%
(\beta,z)\mathrm{Bi}\left(  \vartheta^{^{\frac{2}{3}}}\beta\right)  \right]  ,
\end{align*}
which can be written as%
\begin{gather}
K_{n}(x)\sim\varepsilon^{-\frac{1}{6}}\sqrt{2\pi}\exp\left[  \frac{\left(
z-1\right)  \ln\left(  \frac{U_{0}}{1-z}\right)  +Y^{-}\ln\left(
q-U_{0}\right)  +\left(  1-Y^{-}\right)  \ln\left(  p+U_{0}\right)
+z\ln\left(  \frac{z}{p}\right)  }{\varepsilon}\right] \label{K9}\\
\times\exp\left[  \ln\left(  \frac{U_{0}+p}{q-U_{0}}\right)  \beta
\varepsilon^{-\frac{1}{3}}\right]  \left[  \frac{z\left(  1-z\right)  }%
{pq}\right]  ^{\frac{1}{4}}\vartheta^{-\frac{1}{3}}\left[  \cos\left(  \pi
x\right)  \mathrm{Ai}\left(  \vartheta^{^{\frac{2}{3}}}\beta\right)
-\mathrm{\sin}(\pi x)\mathrm{Bi}\left(  \vartheta^{^{\frac{2}{3}}}%
\beta\right)  \right]  .\nonumber
\end{gather}
Using (\ref{limfac}) and (\ref{limbeta}) in (\ref{K9}) with $\vartheta
=-\Theta,$ we have%
\begin{align}
C_{n}(x)  &  \sim F_{9}(x)=\sqrt{2\pi}\left(  \frac{n}{a}\right)  ^{\frac
{1}{6}}\left(  \sqrt{n}-\sqrt{a}\right)  ^{\frac{1}{3}}\exp\left[  \frac{1}%
{2}n\ln\left(  \frac{n}{a}\right)  +x\ln\left(  \sqrt{\frac{n}{a}}-1\right)
+\sqrt{an}-n\right] \label{C9}\\
&  \times\left\{  \cos\left(  \pi x\right)  \mathrm{Ai}\left[  \left(
\frac{n}{a}\right)  ^{\frac{1}{6}}\frac{\left(  X^{-}-x\right)  }{\left(
\sqrt{n}-\sqrt{a}\right)  ^{\frac{2}{3}}}\right]  -\mathrm{\sin}(\pi
x)\mathrm{Bi}\left[  \left(  \frac{n}{a}\right)  ^{\frac{1}{6}}\frac{\left(
X^{-}-x\right)  }{\left(  \sqrt{n}-\sqrt{a}\right)  ^{\frac{2}{3}}}\right]
\right\} \nonumber
\end{align}
for $x\approx X^{-},\ n>a.$

\item Region X

From (\ref{k10}), (\ref{K3}) and (\ref{K4}) we have for $Y^{-}(z)<y<Y^{+}%
(z),\ 0<z<1$%
\begin{equation}
K_{n}(x)\sim K^{(3)}(y,z)+K^{(4)}(y,z). \label{K10}%
\end{equation}
Thus,%
\begin{equation}
C_{n}(x)\sim F_{10}(x)=F_{3}(x)+F_{4}(x),\quad X^{-}<x<X^{+}. \label{C10}%
\end{equation}

Figure \ref{RegionX} shows the accuracy of the approximation (\ref{C10}) with
$n=30$, $a=2.165184$ in the range $X^{-}<x<X^{+}$.

\begin{figure}[ptb]
\begin{center}
\rotatebox{270} {\resizebox{12cm}{!}{\includegraphics{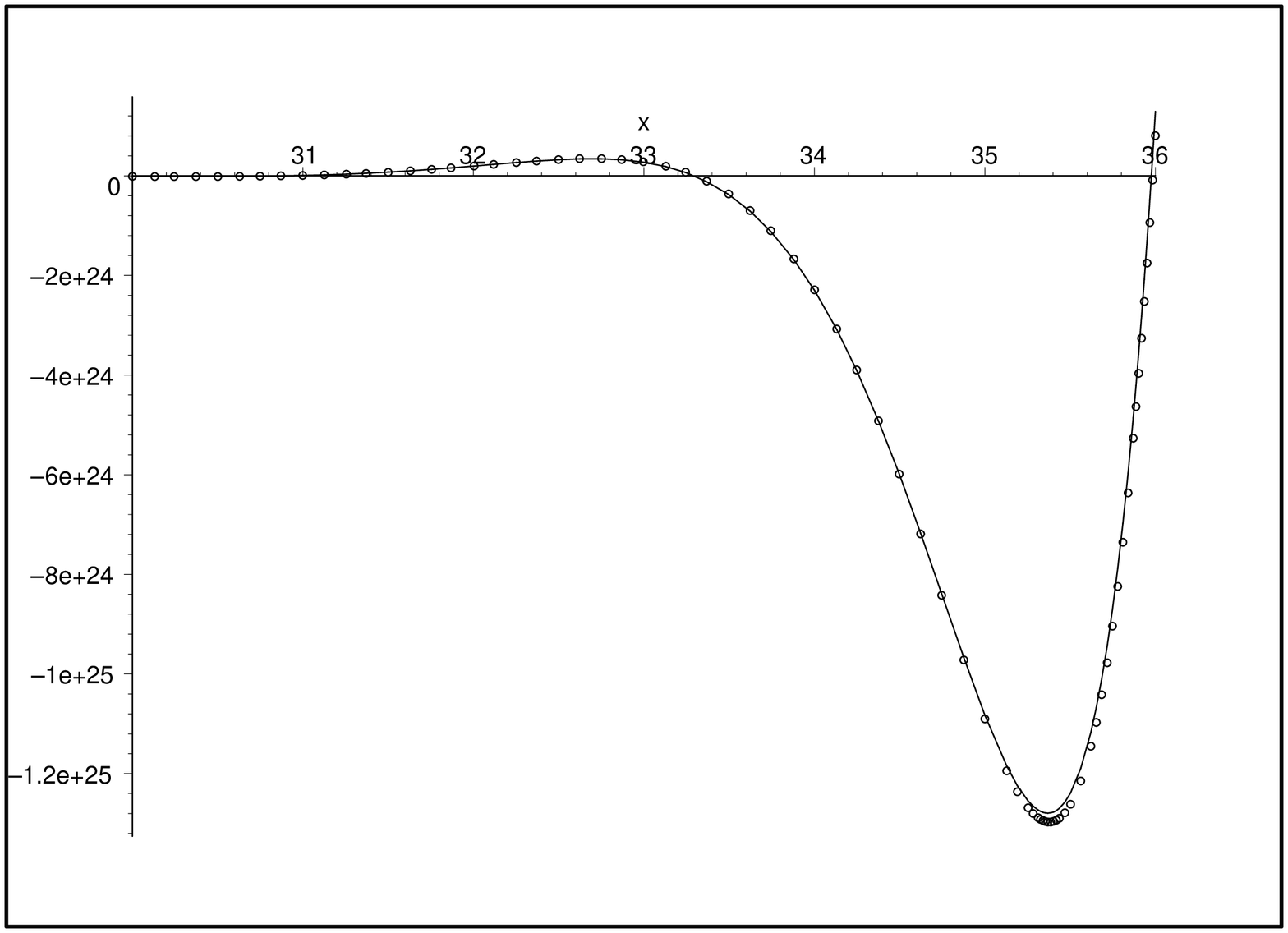}}}
\end{center}
\caption{A comparison of $C_{n}(x)$ (solid curve) and $F_{10}(x)$ (ooo) for $n=30$ with $a=2.165184$.}%
\label{RegionX}%
\end{figure}

\item Region XI

From (\ref{k11}), (\ref{Bi2}) and (\ref{phi}) we have for $y\approx Y^{+}(z),$
$0<z<q,\quad y=Y^{+}(z)+\alpha\varepsilon^{\frac{2}{3}},\quad\alpha=O(1),$%
\begin{align}
K_{n}(x)  &  \sim\varepsilon^{-\frac{1}{6}}\exp\left[  \frac{\psi_{1}%
(z)-\phi(z)}{\varepsilon}+\ln\left(  \frac{U_{0}+q}{U_{0}-p}\right)
\alpha\varepsilon^{-\frac{1}{3}}\right] \label{K11}\\
&  \times\sqrt{2\pi}\left[  \frac{z\left(  1-z\right)  }{pq}\right]
^{\frac{1}{4}}\mathrm{Ai}\left[  \left(  \Theta_{1}\right)  ^{\frac{2}{3}%
}\beta\right]  \left(  \Theta_{1}\right)  ^{-\frac{1}{3}}.\nonumber
\end{align}
Since
\begin{equation}
\alpha=\left[  y-Y^{+}(z)\right]  \varepsilon^{-\frac{2}{3}} \label{alpha}%
\end{equation}
we have from (\ref{Xpm})%
\begin{equation}
\alpha\varepsilon^{-\frac{1}{3}}\rightarrow x-X^{+}. \label{alphalim}%
\end{equation}
From (\ref{psi1}), (\ref{phi}) and (\ref{alphalim}) we get%
\[
\frac{\psi_{1}(z)-\phi(z)}{\varepsilon}+\ln\left(  \frac{U_{0}+q}{U_{0}%
-p}\right)  \alpha\varepsilon^{-\frac{1}{3}}\rightarrow\frac{1}{2}n\ln\left(
\frac{n}{a}\right)  +x\ln\left(  1+\sqrt{\frac{n}{a}}\right)  -\sqrt{an}%
-\sqrt{n}-n\pi\mathrm{i.}%
\]
From (\ref{Theta1}) and (\ref{alphalim}) we obtain%
\begin{equation}
\varepsilon^{-\frac{1}{6}}\sqrt{2\pi}\left[  \frac{z\left(  1-z\right)  }%
{pq}\right]  ^{\frac{1}{4}}\left(  \Theta_{1}\right)  ^{-\frac{1}{3}%
}\rightarrow\sqrt{2\pi}\left(  \frac{n}{a}\right)  ^{\frac{1}{6}}\left(
\sqrt{a}+\sqrt{n}\right)  ^{\frac{1}{3}} \label{limfac1}%
\end{equation}
and%
\begin{equation}
\left(  \Theta_{1}\right)  ^{\frac{2}{3}}\alpha\rightarrow\left(  \frac{n}%
{a}\right)  ^{\frac{1}{6}}\frac{\left(  x-X^{+}\right)  }{\left(  \sqrt
{a}+\sqrt{n}\right)  ^{\frac{2}{3}}}. \label{limalpha}%
\end{equation}
Therefore,%
\begin{align}
C_{n}(x)  &  \sim F_{11}(x)=\sqrt{2\pi}\left(  \frac{n}{a}\right)  ^{\frac
{1}{6}}\left(  \sqrt{a}+\sqrt{n}\right)  ^{\frac{1}{3}}\mathrm{Ai}\left[
\left(  \frac{n}{a}\right)  ^{\frac{1}{6}}\frac{\left(  x-X^{+}\right)
}{\left(  \sqrt{a}+\sqrt{n}\right)  ^{\frac{2}{3}}}\right] \label{C11}\\
&  \times\left(  -1\right)  ^{n}\exp\left[  \frac{1}{2}n\ln\left(  \frac{n}%
{a}\right)  +x\ln\left(  1+\sqrt{\frac{n}{a}}\right)  -\sqrt{an}-\sqrt
{n}\right] \nonumber
\end{align}
for $x\approx X^{+}.$
\end{enumerate}

\section{Comparison with previous results}

We shall now compare our results with those obtained previously in
\cite{MR1297273} and \cite{MR1606887}.

\begin{enumerate}
\item Region VII: $0\leq x<X^{-},$ $n>a.$

Setting $x=un,$ with
\[
u=O(1),\quad0\leq u<1-2\sqrt{\frac{a}{n}}+\frac{a}{n}<1
\]
in (\ref{C7}) we have, as $n\rightarrow\infty$%
\begin{gather}
F_{7}(x)\sim g_{7}(u)=\frac{\cos(un\pi)}{\sqrt{1-u}}\exp\left\{  \left[
u\ln\left(  \frac{n}{a}\right)  +\left(  u-1\right)  \ln\left(  1-u\right)
-u\right]  n+\frac{au}{u-1}\right\} \label{g7}\\
-2\sin\left(  un\pi\right)  \sqrt{\frac{u}{1-u}}\exp\left\{  \left[
\ln\left(  \frac{n}{a}\right)  +\left(  1-u\right)  \ln\left(  1-u\right)
+u\ln(u)-1\right]  n+\frac{a}{1-u}\right\}  .\nonumber
\end{gather}
The second term of equation (\ref{g7}) is the same as the equation before
(5.3) in \cite{MR1297273} and equation (84) in \cite{MR1606887}. However, the
first term is absent in previous works, although it is necessary in the
asymptotic approximation, especially when $u\simeq0,1,2,\ldots.$

\item Region IX: $x\approx X^{-},\ n>a.$

We now set $x=X^{-}+tn^{\frac{1}{6}},\ t=O(1)$ in (\ref{C9}) and obtain, as
$n\rightarrow\infty$%
\begin{gather}
F_{9}(x)\sim g_{9}(t)=\sqrt{2\pi}a^{-\frac{1}{6}}n^{\frac{1}{3}}\exp\left[
\frac{1}{2}\left(  X^{-}+tn^{\frac{1}{6}}+n\right)  \ln\left(  \frac{n}%
{a}\right)  -n+\frac{3}{2}a\right] \label{g9}\\
\times\left\{  \cos\left[  \left(  X^{-}+tn^{\frac{1}{6}}\right)  \pi\right]
\mathrm{Ai}\left(  -ta^{-\frac{1}{6}}\right)  -\sin\left[  \left(
X^{-}+tn^{\frac{1}{6}}\right)  \pi\right]  \mathrm{Bi}\left(  -ta^{-\frac
{1}{6}}\right)  \right\}  .\nonumber
\end{gather}
Equation (\ref{g9}) agrees with equation (5.13) in \cite{MR1297273} and
equation (51) in \cite{MR1606887}.

\item Region X: $X^{-}<x<X^{+}.$

Setting $x=n+a+2\sin(\theta)\sqrt{an},$ with $-\frac{\pi}{2}<\theta<\frac{\pi
}{2}$ in (\ref{C3}) we have, as $n\rightarrow\infty$%
\begin{gather}
F_{3}(x)\sim g_{3}(\theta)=\left(  -1\right)  ^{n}\frac{a^{-\frac{1}{4}%
}n^{\frac{1}{4}}}{\sqrt{2\cos\left(  \theta\right)  }}\exp\left\{  \left[
\ln\left(  \frac{n}{a}\right)  -1\right]  n+\frac{\pi}{4}\mathrm{i}\right\}
\nonumber\\
\times\exp\left\{  \sqrt{an}\left[  \sin\left(  \theta\right)  \ln\left(
\frac{n}{a}\right)  -\sin\left(  \theta\right)  \left(  2\theta-\pi\right)
\mathrm{i}-2\cos\left(  \theta\right)  \mathrm{i}\right]  \right\}
\label{g3}\\
\times\exp\left\{  a\left[  1-\frac{1}{2}\cos\left(  2\theta\right)  +\frac
{1}{2}\ln\left(  \frac{n}{a}\right)  -\frac{1}{2}\sin\left(  2\theta\right)
\mathrm{i}-\theta\mathrm{i}+\frac{\pi}{2}\mathrm{i}\right]  \right\}
.\nonumber
\end{gather}
Similarly from (\ref{C4}) we get%
\begin{gather}
F_{4}(x)\sim g_{4}(\theta)=\left(  -1\right)  ^{n}\frac{a^{-\frac{1}{4}%
}n^{\frac{1}{4}}}{\sqrt{2\cos\left(  \theta\right)  }}\exp\left\{  \left[
\ln\left(  \frac{n}{a}\right)  -1\right]  n-\frac{\pi}{4}\mathrm{i}\right\}
\nonumber\\
\times\exp\left\{  \sqrt{an}\left[  \sin\left(  \theta\right)  \ln\left(
\frac{n}{a}\right)  +\sin\left(  \theta\right)  \left(  2\theta-\pi\right)
\mathrm{i}+2\cos\left(  \theta\right)  \mathrm{i}\right]  \right\}
\label{g4}\\
\times\exp\left\{  a\left[  1-\frac{1}{2}\cos\left(  2\theta\right)  +\frac
{1}{2}\ln\left(  \frac{n}{a}\right)  +\frac{1}{2}\sin\left(  2\theta\right)
\mathrm{i}+\theta\mathrm{i}-\frac{\pi}{2}\mathrm{i}\right]  \right\}
.\nonumber
\end{gather}
Using (\ref{g3}) and (\ref{g4}) in (\ref{C10}) we have%
\begin{gather}
F_{10}(x)\sim g_{10}(\theta)=\left(  -1\right)  ^{n}\frac{\sqrt{2}a^{-\frac
{1}{4}}n^{\frac{1}{4}}}{\sqrt{\cos\left(  \theta\right)  }}\exp\left\{
\left[  \ln\left(  \frac{n}{a}\right)  -1\right]  n\right\} \nonumber\\
\times\exp\left\{  \sqrt{an}\left[  \sin\left(  \theta\right)  \ln\left(
\frac{n}{a}\right)  \right]  +a\left[  1-\frac{1}{2}\cos\left(  2\theta
\right)  +\frac{1}{2}\ln\left(  \frac{n}{a}\right)  \right]  \right\}
\label{g10}\\
\times\cos\left\{  \sqrt{an}\left[  \sin\left(  \theta\right)  \left(
2\theta-\pi\right)  +2\cos\left(  \theta\right)  \right]  +a\left[  \frac
{1}{2}\sin\left(  2\theta\right)  +\theta-\frac{\pi}{2}\right]  -\frac{\pi}%
{4}\right\}  .\nonumber
\end{gather}
Equation (\ref{g10}) is equivalent to equation (44) in \cite{MR1606887}.

\item Region XI: $x\approx X^{+}.$

We now set $x=X^{+}+sn^{\frac{1}{6}},\ s=O(1)$ in (\ref{C11}) and obtain, as
$n\rightarrow\infty$%
\begin{equation}
F_{11}(x)\sim g_{11}(s)=\sqrt{2\pi}a^{-\frac{1}{6}}n^{\frac{1}{3}}\exp\left[
\frac{1}{2}\left(  X^{+}+sn^{\frac{1}{6}}+n\right)  \ln\left(  \frac{n}%
{a}\right)  -n+\frac{3}{2}a\right]  \mathrm{Ai}\left(  sa^{-\frac{1}{6}%
}\right)  . \label{g11}%
\end{equation}
Equation (\ref{g11}) is equation (5.12) in \cite{MR1297273} and equation (30)
in \cite{MR1606887}.
\end{enumerate}

\section{Zeros}

Using the formulas from the previous sections we can obtain approximations to
the zeros of the Charlier polynomials.

\begin{enumerate}
\item $x\simeq0,\ n>a.$

The first zero is exponentially small. From (\ref{C5}) we have, as
$x\rightarrow0$%
\[
C_{5}(x)\sim1+\left[  \ln\left(  \frac{n}{a}-1\right)  -\frac{\sqrt{2\pi n}%
}{n-a}a^{-n}n^{n}e^{a-n}\right]  x.
\]
Solving for $x$ we obtain%
\begin{equation}
x_{0}\simeq\left[  \ln\left(  \frac{n}{a}-1\right)  -\frac{\sqrt{2\pi n}}%
{n-a}a^{-n}n^{n}e^{a-n}\right]  ^{-1}\sim\frac{e^{n-a}a^{n}n^{-n}}{\sqrt{2\pi
n}},\quad n\rightarrow\infty\label{x0}%
\end{equation}
where $x_{0}$ denotes the smallest zero.

\item $0<x<X^{-},\ n>a.$

In this range of $x,$ the zeros are exponentially close to $1,2,\ldots
,\left\lfloor X^{-}\right\rfloor .$ Using
\[
t=\frac{x-n-a+2\sqrt{an}}{n^{\frac{1}{6}}}%
\]
and the asymptotic formulas \cite{atlas}%
\begin{align*}
\mathrm{Ai}(x)  &  \sim\frac{\exp\left[  -\frac{2}{3}x^{\frac{3}{2}}\right]
}{2\sqrt{\pi}x^{\frac{1}{4}}},\quad x\rightarrow\infty\\
\mathrm{Bi}(x)  &  \sim\frac{\exp\left[  \frac{2}{3}x^{\frac{3}{2}}\right]
}{\sqrt{\pi}x^{\frac{1}{4}}},\quad x\rightarrow\infty
\end{align*}
we have, as $n\rightarrow\infty$%
\begin{equation}
\frac{\mathrm{Ai}\left(  -ta^{\frac{1}{6}}\right)  }{\mathrm{Bi}\left(
-ta^{\frac{1}{6}}\right)  }\sim\frac{1}{2}\exp\left[  -\frac{4}{3}a^{-\frac
{1}{4}}n^{-\frac{1}{4}}\left(  X^{-}-x\right)  ^{\frac{3}{2}}\right]  .
\label{AiBi}%
\end{equation}
Using (\ref{AiBi}) in (\ref{g9}) we have%
\[
g_{9}(t)\simeq0\Leftrightarrow\frac{1}{2}\exp\left[  -\frac{4}{3}a^{-\frac
{1}{4}}n^{-\frac{1}{4}}\left(  X^{-}-x\right)  ^{\frac{3}{2}}\right]
\simeq\tan\left(  \pi x\right)  .
\]
Since $x_{j}\simeq j,\ j=1,2,\ldots,\left\lfloor X^{-}\right\rfloor ,$ we get%
\[
\frac{1}{2}\exp\left[  -\frac{4}{3}a^{-\frac{1}{4}}n^{-\frac{1}{4}}\left(
X^{-}-j\right)  ^{\frac{3}{2}}\right]  \simeq\pi\left(  x_{j}-j\right)
\]
which we can solve to obtain%
\begin{equation}
x_{j}\simeq j+\frac{\pi}{2}\exp\left[  -\frac{4}{3}a^{-\frac{1}{4}}%
n^{-\frac{1}{4}}\left(  X^{-}-j\right)  ^{\frac{3}{2}}\right]  ,\quad
\ j=1,2,\ldots,\left\lfloor X^{-}\right\rfloor . \label{xj}%
\end{equation}

\item $X^{-}<x<X^{+}.$

Finally, the non-trivial zeros of the Charlier polynomials can be approximated
using (\ref{g10}). We have $g_{10}(\theta)=0$ \ if and only if%
\[
\cos\left\{  \sqrt{an}\left[  \sin\left(  \theta\right)  \left(  2\theta
-\pi\right)  +2\cos\left(  \theta\right)  \right]  +a\left[  \frac{1}{2}%
\sin\left(  2\theta\right)  +\theta-\frac{\pi}{2}\right]  -\frac{\pi}%
{4}\right\}  =0
\]
or equivalently if%
\[
\sqrt{an}\left[  \sin\left(  \theta\right)  \left(  2\theta-\pi\right)
+2\cos\left(  \theta\right)  \right]  +a\left[  \frac{1}{2}\sin\left(
2\theta\right)  +\theta-\frac{\pi}{2}\right]  -\frac{\pi}{4}=\frac{\pi}{2}+\pi
l,\ l\in\mathbb{Z}%
\]
or%
\begin{equation}
\sqrt{an}\left[  \sin\left(  \theta\right)  \left(  2\theta-\pi\right)
+2\cos\left(  \theta\right)  \right]  +a\left[  \frac{1}{2}\sin\left(
2\theta\right)  +\theta-\frac{\pi}{2}\right]  -\frac{3\pi}{4}-\pi l=0,
\label{equ}%
\end{equation}
with $-\frac{\pi}{2}<\theta<\frac{\pi}{2}.$ Recalling that
\begin{equation}
x=n+a+2\sin(\theta)\sqrt{an}, \label{xtheta}%
\end{equation}
we see that the condition $X^{-}<x<X^{+}$ implies
\begin{equation}
0\leq l\leq2\sqrt{an}-a-\frac{3}{4}. \label{range l}%
\end{equation}
Equation (\ref{equ}) cannot be solved exactly. However, it can be easily
solved numerically to any desired accuracy and using (\ref{xtheta}) gives very
good approximations for the nontrivial zeros.
\end{enumerate}

In Table 1 we computed the exact and approximate zeros of $C_{25}(x)$ with
$a=2.16564899$ using (\ref{x0}), (\ref{xj}) and (\ref{equ})-(\ref{range l}).

\begin{table}[ptb]
\caption{Comparison of the exact and approximate zeros of $C_{25}(x)$ with
$a=2.16564899.$}
\begin{center}
\medskip%
\begin{tabular}
[c]{|c|c|c|}\hline
$l$ & $x$ (exact) & $x$ (approximate)\\\hline
$-$ & $0.41229323\times10^{-16}$ & $0.41549221\times10^{-16}$\\\hline
$-$ & $1.0000000$ & $1.0000000$\\\hline
$-$ & $2.0000000$ & $2.0000000$\\\hline
$-$ & $3.0000000$ & $3.0000001$\\\hline
$-$ & $4.0000000$ & $4.0000009$\\\hline
$-$ & $5.0000001$ & $5.0000073$\\\hline
$-$ & $6.0000015$ & $6.0000507$\\\hline
$-$ & $7.0000227$ & $7.0003063$\\\hline
$-$ & $8.0002574$ & $8.0015785$\\\hline
$-$ & $9.0021153$ & $9.0068260$\\\hline
$-$ & $10.012329$ & $10.024179$\\\hline
$-$ & $11.050278$ & $11.067497$\\\hline
$-$ & $12.147166$ & $12.137242$\\\hline
$11$ & $13.330606$ & $13.334295$\\\hline
$10$ & $14.615276$ & $14.560867$\\\hline
$9$ & $16.007976$ & $15.899727$\\\hline
$8$ & $17.514470$ & $17.350792$\\\hline
$7$ & $19.142918$ & $18.921714$\\\hline
$6$ & $20.905595$ & $20.626110$\\\hline
$5$ & $22.820702$ & $22.484600$\\\hline
$4$ & $24.915443$ & $24.527911$\\\hline
$3$ & $27.232157$ & $26.803591$\\\hline
$2$ & $29.842164$ & $29.391394$\\\hline
$1$ & $32.883964$ & $32.446240$\\\hline
$0$ & $36.717784$ & $36.379078$\\\hline
\end{tabular}
\end{center}
\end{table}

\begin{conclusion}
We analyzed the asymptotic behavior of the Charlier polynomials in the range
$0\leq x$ as $n\rightarrow\infty.$ We also obtained approximations for their
zeros. We intend to extend our method to the other polynomials of the
Askey-scheme to obtain asymptotic expansions of them.
\end{conclusion}

\def\cprime{$'$}

\end{document}